\documentclass[12pt]{amsart}
\usepackage{amsfonts}
\usepackage{eurosym}
\usepackage{amssymb}
\usepackage{amsmath}
\usepackage{enumerate}
\usepackage{graphicx}
\usepackage{caption}
\usepackage{subcaption}
\usepackage{mathrsfs}
\usepackage{float}
\usepackage{lmodern}

\makeatletter
\@namedef{subjclassname@2010}{  \textup{2010} Mathematics Subject Classification}
\makeatother

\newtheorem*{theorem*}{Theorem}

\theoremstyle{definition}

\newtheorem*{sol*}{Solution}

\numberwithin{equation}{section}

\makeatletter
\newcommand{\mathleft}{\@fleqntrue\@mathmargin0pt}
\newcommand{\mathcenter}{\@fleqnfalse}
\makeatother
\begin{document}
\title[IMO problem generalisation]{Generalisation of an IMO Geometry Problem}
\author[D. Kamber Hamzi\'{c}]{Dina Kamber Hamzi\'{c}}
\address{University of Sarajevo \\
Department of Mathematics and Computer Science \\
Zmaja od Bosne 33-35, 71000 Sarajevo, Bosnia and Herzegovina}
\email{dinakamber@pmf.unsa.ba}
\author[L. N\'{e}meth]{L\'{a}szl\'{o} N\'{e}meth}
\address{University of Sopron \\
Institute of Basic Sciences \\
Bajcsy-Zsilinszky 4, 9400 Sopron, Hungary}
\email{nemeth.laszlo@uni-sopron.hu}
\author[Z. \v{S}abanac]{Zenan \v{S}abanac}
\address{University of Sarajevo \\
Department of Mathematics and Computer Science\\
Zmaja od Bosne 33-35, 71000 Sarajevo, Bosnia and Herzegovina}
\email{zsabanac@pmf.unsa.ba}
\maketitle

\begin{abstract}
In this paper, we generalise an interesting geometry problem from the 1995 edition of the International Mathematical Olympiad (IMO) using analytic geometry tools.
\end{abstract}



\section{Introduction and statement of results}

\noindent According to Mitchelmore \cite{MCM}, generalisations are the
cornerstone of school mathematics, covering various aspects like numerical
generalisation in algebra, spatial generalisation in geometry and
measurement, as well as logical generalisations in diverse contexts. The
process of generalising lies at the heart of mathematical activity, serving
as the fundamental method for constructing new knowledge \cite{ET, ELTM}. In
this paper we will generalise an interesting geometry problem appeared in $%
1995$ edition of the International Mathematical Olympiad (IMO) \cite{IMO}.

\bigskip

\begin{figure}[h!]
\centering
\includegraphics[scale=0.95]{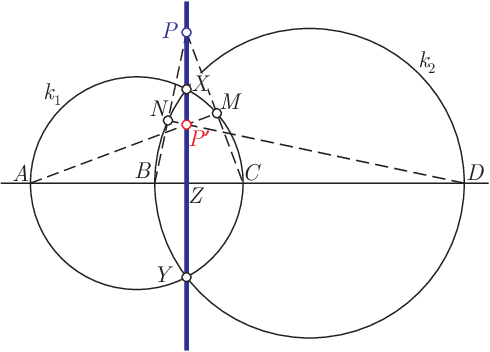}
\caption{36th IMO geometry problem}
\label{fig:IMO_fig_orig}
\end{figure}

\noindent \textbf{Original 36th IMO Geometry Problem.} Let $A,\ B,\ C,\ D$
be four distinct points on a line, in that order. The circles with diameters
$AC $ and $BD$ intersect at $X$ and $Y$. The line $XY$ meets $BC$ at $Z$.
Let $P$ be a point on the line $XY$ other than $Z$. The line $CP$ intersects
the circle with diameter $AC$ at $C$ and $M$, and the line $BP$ intersects
the circle with diameter $BD$ at $B$ and $N$. Prove that the lines $AM,\
DN,\ XY$ are concurrent (Figure~\ref{fig:IMO_fig_orig}).

Its solution, employing analytic geometry, is presented in \cite{KHS}.

\bigskip

Since the line $XY$ in the original problem represents the radical axis of
the two intersecting circles, we shall now expand our investigation to
encompass the radical axis in the broader context. Allow us to reiterate
that the locus of points possessing equal powers concerning two circles is a
line known as the radical axis (or power line) of the circles. Simpler loci
problems should be studied meticulously, as they often pose difficulties
for an average student (see discussion in \cite{Young}). Recent studies have
explored misconceptions that students might have about analytical geometry and
teaching practices in this subject across different educational
levels and contexts (see e.g. \cite{OOK, KFCK, HTST, AAA, CY}).

\bigskip

We will now present a generalised statement of the aforementioned problem
and present its solution using analytical geometry.

\bigskip

\noindent \textbf{Generalisation of the 36th IMO Geometry Problem.} Let $%
k_{1}$ and $k_{2}$ be two circles, and let $\ell $ be the line that contains
their centres. The line $\ell $ and the circle $k_{1}$ intersect at points $%
A $ and $C$, while line $\ell $ and the circle $k_{2}$ intersect at points $%
B $ and $D$. We assume that $A$, $B$, $C$, $D$ appear in that order or in
the order $A$, $C$, $B$, $D$ on $\ell $. Let $p$ be a line perpendicular to
line $\ell $, and let $P$ be any point on line $p$. The line $CP$ intersects
the circle $k_{1}$ at $M$ (in addition to $C$), and the line $BP$ intersects
the circle $k_{2}$ at $N$ (in addition to $B$). The lines $AM$ and $DN$
intersect at $P^{\prime }$. Changing the position of the point $P$ changes
the position of the point $P^{\prime }$. However, $P^{\prime }$ will always
belong to a line $p^{\prime }$, parallel to the line $p$ (Figure~\ref%
{fig:IMO_fig_gen01}). If $p$ is the radical axis of circles $k_{1}$ and $%
k_{2}$, then $p^{\prime }=p$.

\begin{figure}[h!]
\centering
\includegraphics[scale=0.95]{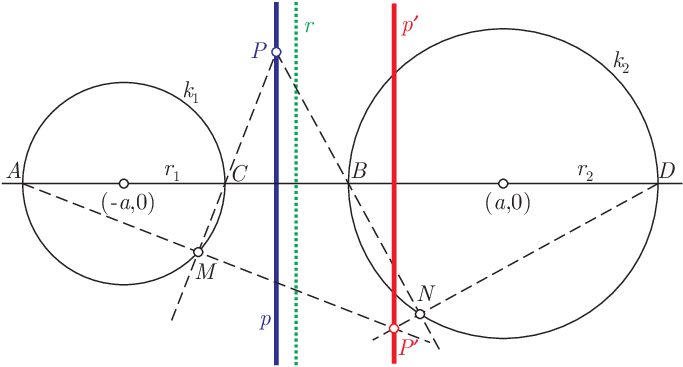}
\caption{36th IMO geometry problem generalisation}
\label{fig:IMO_fig_gen01}
\end{figure}

In the Figure~\ref{fig:IMO_fig_gen01} and all subsequent figures, dotted
line $r$ represents radical axis of two circles.

\section{Analytical Geometry Proof of the Generalisation}

\noindent Without loss of generality, we establish a coordinate system where
$\ell $ represents the $x$-axis, and assume circles $k_{1}$ and $k_{2}$ have
centers symmetrically positioned about the $y$-axis, with $k_{1}$ centered
at $(-a,0)$ with radius $r_{1}$, and $k_{2}$ centered at $(a,0)$ with radius
$r_{2}$. The equations of circles $k_{1}$ and $k_{2}$ are given by ${\left(
x+a\right) }^{2}+y^{2}=r_{1}^{2}$ and ${\left( x-a\right) }%
^{2}+y^{2}=r_{2}^{2}$, respectively. Consequently, we have $A=\left(
-a-r_{1},0\right) $, $B=\left( a-r_{2},0\right) $, $C=\left(
-a+r_{1},0\right) $ and $D=\left( a+r_{2},0\right) $.

\noindent As the line $p$ is perpendicular to the $x$-axis, it is parallel
to the $y$-axis, and we can express it using the equation $x=p$. Henceforth,
the symbol $p$ will serve a dual purpose, representing both the line $p$ and
the $x$-coordinate of point $P$. Consequently, point $P$ can be described by
$P=\left( p,q\right) $, with $q$ changing as point $P$ moves along line $p$.
Let $P^{\prime }=\left( p^{\prime },q^{\prime }\right) $. Using the
tangent-half-parametrisation for the two circles (see, e.g., \cite[para.
5.13 on p. 91]{Rob}), let $M=\left( \frac{r_{1}\left( 1-u^{2}\right) }{%
1+u^{2}}-a,\frac{2r_{1}u}{1+u^{2}}\right) $ and $N=\left( \frac{r_{2}\left(
1-v^{2}\right) }{1+v^{2}}+a,\frac{2r_{2}v}{1+v^{2}}\right) $, for some $%
u,v\in \mathbb{R} \cup \left\{\infty\right\} $.

\bigskip

\noindent Now, we examine the condition for points $P$, $C$ and $M$ to be
collinear by considering the vanishing of the determinant (\cite[para. 1.72
on pp. 14-15]{Rob})
\begin{equation}
\begin{vmatrix}
p & q & 1 \\
-a+r_{1} & 0 & 1 \\
\frac{r_{1}\left( 1-u^{2}\right) }{1+u^{2}}-a & \frac{2r_{1}u}{1+u^{2}} & 1%
\end{vmatrix}%
.  \label{det1}
\end{equation}%
Subtracting the second row from the third row, and then taking out the
common factors from the third row, we have%
\begin{equation*}
\frac{2r_{1}u}{1+u^{2}}%
\begin{vmatrix}
p & q & 1 \\
-a+r_{1} & 0 & 1 \\
-u & 1 & 0%
\end{vmatrix}%
=0.
\end{equation*}%
Expanding the determinant along the third row, we obtain%
\begin{equation}
\frac{2r_{1}u}{1+u^{2}}\left( -uq-p-a+r_{1}\right) =0.  \label{eq_u1}
\end{equation}%
Solving for $u$, and discarding the solution $u=0$, which would entail $M=C$%
, we get%
\begin{equation}
u=\frac{r_{1}-a-p}{q}.  \label{eq_u2}
\end{equation}%
The cases $u=0$ and $q=0$, where the latter leads to $u=\infty $
and $M=A$, will be addressed later, along with the cases when the
denominator becomes zero or tends to infinity.

\noindent Similarly, the condition for the collinearity of $P^{\prime }$, $A$
and $M$ is the vanishing of the determinant%
\begin{equation*}
\begin{vmatrix}
p^{\prime } & q^{\prime } & 1 \\
-a-r_{1} & 0 & 1 \\
\frac{r_{1}\left( 1-u^{2}\right) }{1+u^{2}}-a & \frac{2r_{1}u}{1+u^{2}} & 1%
\end{vmatrix}%
.
\end{equation*}%
The last determinant can be obtained from \eqref{det1} by replacing $p$, $q$%
, $r_{1}$ and $u$ by $p^{\prime }$, $q^{\prime }$, $-r_{1}$ and $-u^{-1}$,
respectively. Applying this to \eqref{eq_u1} leads to
\begin{equation*}
\frac{2r_{1}}{1+u^{2}}\left( q^{\prime }-up^{\prime }-ua-ur_{1}\right) =0.
\end{equation*}%
Solving for $u$, and discarding the solution $u=\infty $, which would entail
$M=A$, we get%
\begin{equation*}
u=\frac{q^{\prime }}{r_{1}+a+p^{\prime }}.
\end{equation*}%
From this and \eqref{eq_u2} we now have%
\begin{equation}
\frac{r_{1}-a-p}{q}=\frac{q^{\prime }}{r_{1}+a+p^{\prime }}.  \label{eq1}
\end{equation}%
Likewise, the conditions for collinearity among points $P$, $B$, $N$, and
of $P^{\prime }$, $D$, $N$, are determined by the vanishing of the
respective determinants%
\begin{equation*}
\begin{vmatrix}
p & q & 1 \\
a-r_{2} & 0 & 1 \\
\frac{r_{2}\left( 1-v^{2}\right) }{1+v^{2}}+a & \frac{2r_{2}v}{1+v^{2}} & 1%
\end{vmatrix}%
\text{ and }%
\begin{vmatrix}
p^{\prime } & q^{\prime } & 1 \\
a+r_{2} & 0 & 1 \\
\frac{r_{2}\left( 1-v^{2}\right) }{1+v^{2}}+a & \frac{2r_{2}v}{1+v^{2}} & 1%
\end{vmatrix}%
\text{.}
\end{equation*}%
The first one can be obtained from \eqref{det1} by replacing $a$, $r_{1}$
and $u$ by $-a$, $-r_{2}$ and $-v^{-1}$, respectively; and the second one by
replacing $p$, $q$, $a$, $r_{1}$ and $u$ by $p^{\prime }$, $q^{\prime }$, $%
-a $, $r_{2}$ and $v$, respectively. Applying the same changes to %
\eqref{eq_u1}, we obtain $v=\frac{q}{r_{2}-a+p}$ and $v=\frac{%
r_{2}+a-p^{\prime }}{q^{\prime }}$, respectively. Hence,%
\begin{equation}
\frac{q}{r_{2}-a+p}=\frac{r_{2}+a-p^{\prime }}{q^{\prime }}.  \label{eq2}
\end{equation}%
Multiplying \eqref{eq1} and \eqref{eq2} yields the following relationship
\begin{equation}
\frac{r_{1}-a-p}{r_{2}-a+p}=\frac{r_{2}+a-p^{\prime }}{r_{1}+a+p^{\prime }},
\label{eq3}
\end{equation}%
which indicates that $p^{\prime }$ depends on $p$ (as well as on $a$, $r_{1}$
and $r_{2} $) but is independent of $q$. This results in $P^{\prime}$ lying
on a line perpendicular to the $x$-axis and, thus, parallel to line $p$,
which was what we wanted to prove.

To get an explicit value for $p^{\prime }$, we perform crosswise
multiplication in \eqref{eq3}, and obtain%
\begin{equation*}
\left( r_{1}-a-p\right) \left( r_{1}+a+p^{\prime }\right) =\left(
r_{2}-a+p\right) \left( r_{2}+a-p^{\prime }\right) .
\end{equation*}%
This yields%
\begin{gather*}
\left( r_{1}-a-p\right) p^{\prime }+\left( r_{2}-a+p\right) p^{\prime
}=\left( r_{2}-a+p\right) \left( r_{2}+a\right) -\left( r_{1}-a-p\right)
\left( r_{1}+a\right) \\
=\left( r_{2}-a\right) \left( r_{2}+a\right) +p\left( r_{2}+a\right) -\left(
r_{1}-a\right) \left( r_{1}+a\right) +p\left( r_{1}+a\right) .
\end{gather*}%
From this we obtain%
\begin{equation*}
\left( r_{1}+r_{2}-2a\right) p^{\prime }=r_{2}^{2}-r_{1}^{2}+p\left(
r_{1}+r_{2}+2a\right) ,
\end{equation*}%
and so%
\begin{equation}
p^{\prime }=\frac{r_{2}^{2}-r_{1}^{2}+p\left( r_{1}+r_{2}+2a\right) }{%
r_{1}+r_{2}-2a}.  \label{eq4}
\end{equation}%
To get $q{\prime }$, we take the reciprocal of \eqref{eq1} to obtain $\frac{q}{r_{1}-a-p}=%
\frac{r_{1}+a+p^{\prime }}{q^{\prime }}$, and upon adding this to \eqref{eq2}%
, we arrive at%
\begin{equation*}
q\left( \frac{1}{r_{1}-a-p}+\frac{1}{r_{2}-a+p}\right) =\frac{r_{1}+r_{2}+2a%
}{q^{\prime }}.
\end{equation*}%
Now,
\begin{equation}
q^{\prime }=\frac{\left( r_{1}+r_{2}+2a\right) \left( a-r_{1}+p\right)
\left( a-r_{2}-p\right) }{q\left( r_{1}+r_{2}-2a\right) }.  \label{eq5}
\end{equation}%
From \eqref{eq4} and \eqref{eq5}, we see that
\begin{equation*}
P^{\prime }=\left( \frac{r_{2}^{2}-r_{1}^{2}+p\left( r_{1}+r_{2}+2a\right) }{%
r_{1}+r_{2}-2a},\frac{\left( r_{1}+r_{2}+2a\right) \left( a-r_{1}+p\right)
\left( a-r_{2}-p\right) }{q\left( r_{1}+r_{2}-2a\right) }\right) .
\end{equation*}

\bigskip

Now, we will consider special cases.

\bigskip

\begin{figure}[h!]
\centering
\includegraphics[scale=0.90]{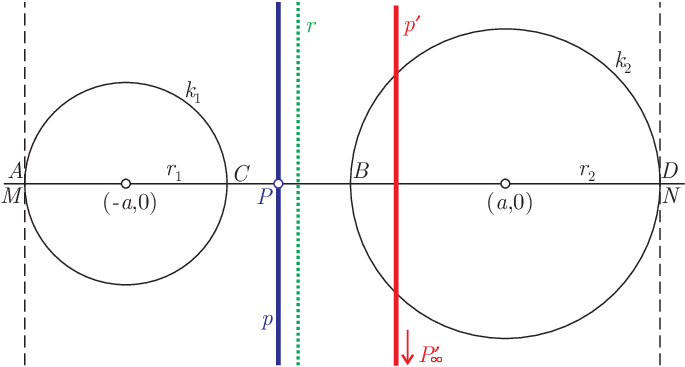}
\caption{$P^{\prime}$ is at infinity}
\label{fig:IMO_fig_gen03}
\end{figure}

If $P$ is on the $x$-axis, i.e., if $q=0$, then $q^{\prime }=\infty $, and $%
P^{\prime }$ is the point at infinity. Then, by the equations $u=\frac{%
r_{1}-a-p}{q}$ and $v=\frac{q}{r_{2}-a+p}$, we have $u=\infty $ and $v=0$,
i.e., $M=A$ and $N=D$. The tangent lines of the circles $k_{1}$ and $k_{2}$
at points $A$ and $D$, respectively, represent the limiting position of the
lines $AM$ and $DN$ in the previous case (see Figure~\ref{fig:IMO_fig_gen03}%
). Their intersection point disappears (or we can say that they intersect at
the infinite point $P_{\infty }^{\prime }$ on the line $p^{\prime }$).

\begin{figure}[h!]
\centering
\includegraphics[scale=0.85]{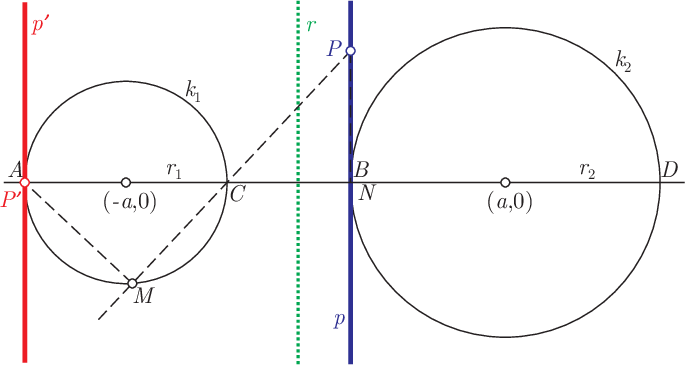}
\caption{$p$ touching $k_{2}$ at point $B$}
\label{fig:IMO_fig_gen02}
\end{figure}

If the line $x=p$ passes through the points $B$ or $C$, then $p=a-r_{2}$ or $%
r_{1}-a$, respectively, and by \eqref{eq3}, $p^{\prime }=-a-r_{1}$ or $%
a+r_{2}$, respectively. Then, by \eqref{eq5}, $q^{\prime }=0$ in either
case, so $P^{\prime }=A$ or $D$, respectively, and loci of point $P^{\prime }
$ consist of a single point (see Figure~\ref{fig:IMO_fig_gen02}). Within
lines passing through the point $A$, resp. $D$, there is obviously the line $%
p^{\prime }$ tangent to the circle $k_{1}$, resp. $k_{2}$, and parallel to $p$.

Next, if $r_{1}+r_{2}=2a$, the circles $k_{1}$ and $k_{2}$ touch externally
(in which case $B=C$), and by \eqref{eq4} and \eqref{eq5} we have $p^{\prime
}=\infty $ and $q^{\prime }=\infty $; therefore, point $P^{\prime }$ is at
infinity. Thus, lines $AM$ and $ND$ are parallel (see Figure~\ref%
{fig:IMO_fig_touch}), resulting in the vanishing of line $p^{\prime }$. Notice that one may conclude the same from the fact that the points $P$, $C$, $B$, $M$, $N$ are collinear, and $\angle {AMC}%
=\angle {DNB}=90^{\circ }$ as angles subtended by the respective
diameters.
\begin{figure}[h!]
\centering
\includegraphics[scale=0.95]{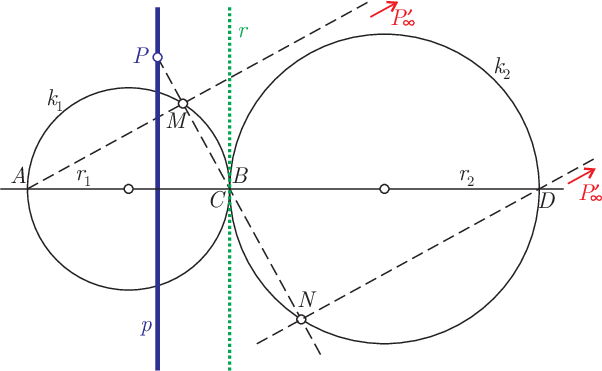}
\caption{Touching circles}
\label{fig:IMO_fig_touch}
\end{figure}

Finally, the equation defining the radical axis for circles $k_{1}$ and $%
k_{2}$ gives us (see \cite[para. 5.51 on p. 99]{Rob})
\begin{equation*}
{\left( x+a\right) }^{2}+y^{2}-r_{1}^{2}={\left( x-a\right) }%
^{2}+y^{2}-r_{2}^{2},
\end{equation*}%
which simplifies to
\begin{equation*}
x=\frac{r_{1}^{2}-r_{2}^{2}}{4a}.
\end{equation*}%
From \eqref{eq4}, $p=p^{\prime }$ if, and only if, $p=\frac{%
r_{1}^{2}-r_{2}^{2}}{4a}$, which means that if $p$ is the radical axis of
the two circles $k_{1}$ and $k_{2}$, then $p^{\prime }=p$.

\bigskip

\noindent \textit{Remark.} In case that points on line $\ell $ appear in
order $A$, $B$, $C$, $D$, and line $p$ passes through the intersection
points of the two circles (in which case line $p$ is the radical axis of the
circles), we obtain the original IMO problem.

\bigskip

Note that more cases can be considered by changing the roles of points $A$ and $C$,
and/or $B$ and $D$, placing one circle inside the other, or having one touch
the other internally. We invite interested reader to investigate these other
cases.

\end{document}